\documentclass[a4paper,twoside,11pt]{article}
\usepackage{amsfonts}
\usepackage{amssymb}\usepackage{amsmath}
\newtheorem{thm}{Theorem}[section]
 \newtheorem{cor}{Corollary}[section]

 \numberwithin{equation}{section}
\newcommand{\double}{\baselineskip 1.24 \baselineskip}

\title{Generalizations of H\"{o}lder's and some related integral inequalities on fractal
space
}

\author{{Guang-Sheng  Chen\thanks{\text{E-mail address}: cgswavelets@126.com(Chen)
}\quad}\\
{\small Department of Computer Engineering, Guangxi Modern
Vocational Technology College,} \\{\small Hechi,Guangxi, 547000,
P.R. China}
}
\begin{document}
\date{}
\maketitle \double

\textbf{Abstract:}\quad Based on the local fractional calculus, we establish some
new generalizations of H\"{o}lder's inequality. By using it, some results on the
generalized integral inequality in fractal space are investigated in detail.   

\section{Introduction}
\hskip\parindent
While the renowned inequality of H\"{o}lder [1] is well celebrated for its
beauty and its wide range of important applications to real and complex
analysis, functional analysis, as well as many disciplines in applied
mathematics. The purpose of this work is to establish some generalizations
of H\"{o}lder inequality on local fractional calculus and other inequality
based on it. Fractal calculus (also called local fractional calculus) has
played an important role in not only mathematics but also in physics and
engineers [2-15]. Local fractional derivative [6-8] were written in the form

$f^{(\alpha )}(x_0 ) = \frac{d^\alpha f(x)}{dx^\alpha }\vert _{x = x_0 } =
\mathop {\lim }\limits_{x \to x_0 } \frac{\Delta ^\alpha (f(x) - f(x_0
))}{(x - x_0 )^\alpha }$ for $0 < \alpha \le 1$,\\
where
\[
\Delta ^\alpha (f(x) - f(x_0 )) \cong \Gamma (1 + \alpha )\Delta (f(x) -
f(x_0 )).
\]
Local fractional integral of $f(x)$ [6-7,9] was written in the form
\[
{ }_aI_b^{(\alpha )} f(x) = \frac{1}{\Gamma (1 + \alpha )}\int_a^b
{f(t)(dt)^\alpha } = \frac{1}{\Gamma (1 + \alpha )}\mathop {\lim
}\limits_{\Delta t \to 0} \sum\limits_{j = 0}^{N - 1} {f(t_j )(\Delta t_j
)^\alpha } ,
\]
\\ with $\Delta t_j = t_{j + 1} - t_j $ and $\Delta t = \max \{\Delta t_1 ,\Delta
t_2 ,\ldots ,\Delta t_j ,\ldots \}$, where for $j = 1,2,\ldots ,N - 1$，$t_0
= a$ and $t_N = b$, $[t_j ,t_{j + 1} ]$ is a partition of the interval $[a,b]$.
Aims of this paper are to study the some new generalizations of H\"{o}lder's
inequality and some results based on them.

\section{Some generalizations of H\"{o}lder inequality}
\hskip\parindent
In the section, we give some generalizations of H\"{o}lder inequality. In order to prove our results, we first review the H\"{o}lder inequality [15]:
\begin{thm}
\label{Theorem 2.1}.[15] Let $f(x)$, $g(x) \ge 0$, $p > 1$, $\frac{1}{p} +
\frac{1}{q} = 1$ , then
\begin{equation*}
\frac{1}{\Gamma (1 + \alpha )}\int_a^b {\left| {f(x)g(x)} \right|(dx)^\alpha
} \le \left( {\frac{1}{\Gamma (1 + \alpha )}\int_a^b {\left| {f(x)}
\right|^p(dx)^\alpha } } \right)^{1 / p}\left( {\frac{1}{\Gamma (1 + \alpha
)}\int_a^b {\left| {g(x)} \right|^q(dx)^\alpha } } \right)^{1 / q},
\tag{2.1}
\label{2.1}
\end{equation*}
equalities holding if and only if $f(x) = \lambda g(x)$, where $\lambda $ is
a constant.
\end{thm}
Based on Theorem 2.1, we have the following important result.
\begin{thm}\label{Theorem 2.2.} Let $f(x)$, $g(x) \ge 0$, $0 < p < 1$, $\frac{1}{p} +
\frac{1}{q} = 1$ , then
\begin{equation*}
\frac{1}{\Gamma (1 + \alpha )}\int_a^b {\left| {f(x)g(x)} \right|(dx)^\alpha
} \ge \left( {\frac{1}{\Gamma (1 + \alpha )}\int_a^b {\left| {f(x)}
\right|^p(dx)^\alpha } } \right)^{1 / p}\left( {\frac{1}{\Gamma (1 + \alpha
)}\int_a^b {\left| {g(x)} \right|^q(dx)^\alpha } } \right)^{1 / q},
\tag{2.2}
\label{2.2}
\end{equation*}
equalities holding if and only if $f(x) = \lambda g(x)$, where $\lambda $ is
a constant.
\end{thm}
\textbf{Proof. }Set $c = 1 / p$, then we have $q = - pd$, $d = c / (c - 1)$.
By (2.1), we obtain
\begin{equation*}
\begin{split}
 &\frac{1}{\Gamma (1 + \alpha )}\int_a^b {\left| {f(x)} \right|^p(dx)^\alpha
} = \frac{1}{\Gamma (1 + \alpha )}\int_a^b {\left| {f(x)g(x)}
\right|^p\left| {g(x)} \right|^{ - p}(dx)^\alpha } \\
 &\le \left( {\frac{1}{\Gamma (1 + \alpha )}\int_a^b {\left| {f(x)g(x)}
\right|^{pc}(dx)^\alpha } } \right)^{1 / c}\left( {\frac{1}{\Gamma (1 +
\alpha )}\int_a^b {\left| {g(x)} \right|^{ - pd}(dx)^\alpha } } \right)^{1 /
d} \\
 &= \left( {\frac{1}{\Gamma (1 + \alpha )}\int_a^b {\left| {f(x)g(x)}
\right|(dx)^\alpha } } \right)^{1 / c}\left( {\frac{1}{\Gamma (1 + \alpha
)}\int_a^b {\left| {g(x)} \right|^q(dx)^\alpha } } \right)^{1 - p}. \\
 \end{split}
\tag{2.3}
\label{2.3}
\end{equation*}
In (2.3), multiplying both sides by $\left( {\frac{1}{\Gamma (1 + \alpha
)}\int_a^b {\left| {g(x)} \right|^q(dx)^\alpha } } \right)^{p - 1}$ yields
\begin{equation*}
\frac{1}{\Gamma (1 + \alpha )}\int_a^b {\left| {f(x)} \right|^p(dx)^\alpha }
\left( {\frac{1}{\Gamma (1 + \alpha )}\int_a^b {\left| {g(x)}
\right|^q(dx)^\alpha } } \right)^{p - 1} \le \left( {\frac{1}{\Gamma (1 +
\alpha )}\int_a^b {\left| {f(x)g(x)} \right|(dx)^\alpha } } \right)^p.
\tag{2.4}
\label{2.4}
\end{equation*}
Using (2.4) implies that
\[
\frac{1}{\Gamma (1 + \alpha )}\int_a^b {\left| {f(x)g(x)} \right|(dx)^\alpha
} \ge \left( {\frac{1}{\Gamma (1 + \alpha )}\int_a^b {\left| {f(x)}
\right|^p(dx)^\alpha } } \right)^{1 / p}\left( {\frac{1}{\Gamma (1 + \alpha
)}\int_a^b {\left| {g(x)} \right|^q(dx)^\alpha } } \right)^{1 / q}.
\]

By Theorem 2.1 and Theorem 2.2, we give the following result.
\begin{cor}
\label{Corollary 2.3.} Let $f_j (x) \ge 0$, $p_j \in {\rm R}$ , $j =
1,2,\ldots m$, $\sum\limits_{j = 1}^m {1 / p_j } = 1$.

(1) for $p_j > 1$, we have
\begin{equation*}
\frac{1}{\Gamma (1 + \alpha )}\int_a^b {\prod\limits_{j = 1}^m {\left| {f_j
(x)} \right|(dx)^\alpha } } \le \prod\limits_{j = 1}^m {(\int_a^b
{\frac{1}{\Gamma (1 + \alpha )}\left| {f_j (x)} \right|^{p_j }(dx)^\alpha }
)^{1 / p_j }} .
\tag{2.5}\label{2.5}
\end{equation*}

(2) for $0 < p_1 < 1$, $p_j < 0$, $j = 2,\ldots m$,we have
\begin{equation*}
\frac{1}{\Gamma (1 + \alpha )}\int_a^b {\prod\limits_{j = 1}^m {\left| {f_j
(x)} \right|(dx)^\alpha } } \ge \prod\limits_{j = 1}^m {(\int_a^b
{\frac{1}{\Gamma (1 + \alpha )}\left| {f_j (x)} \right|^{p_j }(dx)^\alpha }
)^{1 / p_j }}.
\tag{2.6}\label{2.6}
\end{equation*}
\end{cor}
\textbf{Proof.} (1)We use induction on $m$. When $m = 2$, we are given $p_1 $, $p_2
> 0$ with $1 / p_1 + 1 / p_2 = 1$. In particular, we have $p_1 $, $p{ }_2 >
1$ and so (2.5) is reduced to the H\"{o}der's inequality (2.1). Now suppose
(2.5) holds for some integer $m \ge 2$. We claim that it also holds for $m +
1$. So let $p_1 ,p{ }_2,\ldots ,p{ }_{m + 1} > 0$ be real numbers with
$\sum\limits_{j = 1}^{m + 1} {1 / p_j } = 1$ and let $f_j (x) \ge 0$, $j =
1,2,\ldots m,m + 1$. Note that, as above, we must have $p_i > 1$ for $j =
1,2,\ldots m,m + 1$. In particular, we have
\begin{equation*}
p_1 > 0,
p_1 / (p_1 - 1) > 0,
1 / p_1 + (p_1 - 1) / p_1 = 1.
\tag{2.7}
\label{2.7}
\end{equation*}

Thus by the H\"{o}der's inequality (2.1),
\begin{equation*}
\begin{split}
 &\frac{1}{\Gamma (1 + \alpha )}\int_a^b {\prod\limits_{j = 1}^{m + 1}
{\left| {f_j (x)} \right|(dx)^\alpha } } = \frac{1}{\Gamma (1 + \alpha
)}\int_a^b {\left| {f_1 (x)} \right|\prod\limits_{j = 2}^{m + 1} {\left|
{f_j (x)} \right|(dx)^\alpha } } \\
 &\le (\frac{1}{\Gamma (1 + \alpha )}\int_a^b {\left| {f_1 (x)} \right|^{p_1
}(dx)^\alpha } )^{1 / p_1 }\left( {\frac{1}{\Gamma (1 + \alpha )}\int_a^b
{\left( {\prod\limits_{j = 2}^{m + 1} {\left| {f_j (x)} \right|} }
\right)^{p_1 / (p_1 - 1)}(dx)^\alpha } } \right)^{(p_1 - 1) / p_1 } \\
 &= (\frac{1}{\Gamma (1 + \alpha )}\int_a^b {\left| {f_1 (x)} \right|^{p_1
}(dx)^\alpha } )^{1 / p_1 }\left( {\frac{1}{\Gamma (1 + \alpha )}\int_a^b
{\prod\limits_{j = 2}^{m + 1} {\left| {f_j (x)} \right|^{p_1 / (p_1 - 1)}}
(dx)^\alpha } } \right)^{(p_1 - 1) / p_1 }.\\
 \end{split}
 \tag{2.8}
 \label{2.8}
\end{equation*}
Next, since
\begin{equation*}
p_j (p_1 - 1) / p_1 > 0, \text{for} \quad j = 2,\ldots m,m + 1 . \tag{2.9}
 \label{2.9}
\end{equation*}
\begin{equation*}
\sum\limits_{j = 2}^{m + 1} {p_1 / p_j (p_1 - 1)} = p_1 (p_1 -
1)\sum\limits_{j = 2}^{m + 1} {1 / p_j } = p_1 (p_1 - 1)(1 - 1 / p_1 ) = 1.
 \tag{2.10}
 \label{2.10}
\end{equation*}
\noindent
by induction hypothesis and (2.8), we arrive at
\begin{equation*}
\begin{split}
& \frac{1}{\Gamma (1 + \alpha )}\int_a^b {\prod\limits_{j = 1}^{m + 1}
{\left| {f_j (x)} \right|} (dx)^\alpha } \le (\frac{1}{\Gamma (1 + \alpha
)}\int_a^b {\left| {f_1 (x)} \right|^{p_1 }(dx)^\alpha } )^{1 / p_1 } \\
 &\times \left( {\prod\limits_{j = 2}^{m + 1} {\left( {\frac{1}{\Gamma (1 +
\alpha )}\int_a^b {\left| {f_j (x)} \right|^{p_1 / (p_1 - 1) \cdot p_j (p_1
- 1) / p_1 }(dx)^\alpha } } \right)^{p_1 / p_j (p_1 - 1)}} } \right)^{(p_1 -
1) / p_1 } \\
& = (\frac{1}{\Gamma (1 + \alpha )}\int_a^b {\left| {f_1 (x)} \right|^{p_1
}(dx)^\alpha } )^{1 / p_1 }\prod\limits_{j = 2}^{m + 1} {\left(
{\frac{1}{\Gamma (1 + \alpha )}\int_a^b {\left| {f_j (x)} \right|^{p_j
}(dx)^\alpha } } \right)^{1 / p_j }}. \\
\end{split}
 \tag{2.11}
 \label{2.11}
\end{equation*}
Hence, we arrive at the result.

(2) Similar to the proof of (2.5), we use induction on $m$. Clearly
when $m = 2$, equation (2.6) reduces to the H\"{o}lder's inequality (2.2).
Now suppose that (2.6) holds for some integer $m \ge 2$. We claim that it
also holds for $m + 1$. So let

$p_1 ,p{ }_2,\ldots ,p_m < 0$ and $p_{m + 1} \in {\rm R}$ be such that
$\sum\limits_{j = 1}^{m + 1} {1 / p_j } = 1$ and let $f_j (x)\geq 0$, $j = 1,2,\ldots m,m + 1$. Note that $0 < p_{m + 1} <
1$, since
\begin{equation*}
p_1 > 0,
0 < p_1 / (p_1 - 1) < 1,
1 / p_1 + (p_1 - 1) / p_1 = 1.
\tag{2.12}
\label{2.12}
\end{equation*}
by the H\"{o}der's inequality (2.2), we have
\begin{equation*}
\begin{split}
 &\frac{1}{\Gamma (1 + \alpha )}\int_a^b {\prod\limits_{j = 1}^m {\left| {f_j
(x)} \right|} (dx)^\alpha } = \frac{1}{\Gamma (1 + \alpha )}\int_a^b {\left|
{f_1 (x)} \right|\prod\limits_{j = 2}^m {\left| {f_j (x)} \right|}
(dx)^\alpha } \\
 &\ge \left( {\frac{1}{\Gamma (1 + \alpha )}\int_a^b {\left| {f_1 (x)}
\right|^{p_1 }(dx)^\alpha } } \right)^{1 / p_1 }\left( {\frac{1}{\Gamma (1 +
\alpha )}\int_a^b {\left( {\prod\limits_{j = 2}^m {\left| {f_j (x)} \right|}
} \right)^{p_1 / (p_1 - 1)}(dx)^\alpha } } \right)^{(p_1 - 1) / p_1 } \\
 &= (\frac{1}{\Gamma (1 + \alpha )}\int_a^b {\left| {f_1 (x)} \right|^{p_1
}(dx)^\alpha } )^{1 / p_1 }\left( {\frac{1}{\Gamma (1 + \alpha )}\int_a^b
{\prod\limits_{j = 2}^m {\left| {f_j (x)} \right|^{p_1 / (p_1 - 1)}}
(dx)^\alpha } } \right)^{(p_1 - 1) / p_1 }. \\
 \end{split}
 \tag{2.13}
 \label{2.13}
\end{equation*}
unless $\frac{1}{\Gamma (1 + \alpha )}\int_a^b {\left| {f_1 (x)}
\right|^{p_1 }(dx)^\alpha } = 0$.

Now since
\begin{equation*}p_j (p_1 - 1) / p_1 < 0 \text{for} j = 2,\ldots m, p_{m + 1} (p_1 - 1) / p_1 >0
\tag{2.14}
 \label{2.14}
\end{equation*}
and as in (2.10),
\begin{equation*}
\sum\limits_{j = 2}^{m + 1} {p_1 / p_j (p_1 - 1)} = 1.
\tag{2.15}
 \label{2.15}
\end{equation*}
by induction hypothesis and (2.13), we obtain
\begin{equation*}
\begin{split}
 &\frac{1}{\Gamma (1 + \alpha )}\int_a^b {\prod\limits_{j = 1}^{m + 1}
{\left| {f_j (x)} \right|} (dx)^\alpha } \ge \left( {\frac{1}{\Gamma (1 +
\alpha )}\int_a^b {\left| {f_1 (x)} \right|^{p_1 }(dx)^\alpha } } \right)^{1
/ p_1 } \\
 &\times \left( {\prod\limits_{j = 2}^{m + 1} {\left( {\frac{1}{\Gamma (1 +
\alpha )}\int_a^b {\left| {f_j (x)} \right|^{p_1 / (p_1 - 1) \cdot p_j (p_1
- 1) / p_1 }(dx)^\alpha } } \right)^{p_1 / p_j (p_1 - 1)}} } \right)^{(p_1 -
1) / p_1 } \\
 &= \left( {\frac{1}{\Gamma (1 + \alpha )}\int_a^b {\left| {f_1 (x)}
\right|^{p_1 }(dx)^\alpha } } \right)^{1 / p_1 }\prod\limits_{j = 2}^{m + 1}
{\left( {\frac{1}{\Gamma (1 + \alpha )}\int_a^b {\left| {f_j (x)}
\right|^{p_j }(dx)^\alpha } } \right)^{1 / p_j }} \\
 &= \prod\limits_{j = 1}^{m + 1} {\left( {\frac{1}{\Gamma (1 + \alpha
)}\int_a^b {\left| {f_j (x)} \right|^{p_j }(dx)^\alpha } } \right)^{1 / p_j
}} .\\
\end{split}
 \tag{2.16}
 \label{2.16}
\end{equation*}
unless $\frac{1}{\Gamma (1 + \alpha )}\int_a^b {\left| {f_j (x)}
\right|^{p_j }(dx)^\alpha } = 0$ for some $j = 1,2,\ldots m$.

\section{Some related results }
\hskip\parindent
To set the stage, we recall Minkowski inequality [15]:
\begin{thm}\label{Theorem 3.1}.[15] Let $f(x)$, $g(x) \ge 0$, $p > 1$ , then
\begin{equation*}
\begin{split}
 &\left( {\frac{1}{\Gamma (1 + \alpha )}\int_a^b {\left| {f(x) + g(x)}
\right|^p(dx)^\alpha } } \right)^{1 / p} \\
 &\le \left( {\frac{1}{\Gamma (1 + \alpha )}\int_a^b {\left| {f(x)}
\right|^p(dx)^\alpha } } \right)^{1 / p} + \left( {\frac{1}{\Gamma (1 +
\alpha )}\int_a^b {\left| {g(x)} \right|^q(dx)^\alpha } } \right)^{1 / q}. \\
 \end{split}
 \tag{3.1}
 \label{3.1}
\end{equation*}
equalities holding if and only if $f(x) = \lambda g(x)$, where $\lambda $ is
a constant.
\end{thm}
\begin{thm}\label{Theorem 3.2}.Let $f(x)$, $g(x) \ge 0$, $0 < p < 1$ , then
\begin{equation*}
\begin{split}
& \left( {\frac{1}{\Gamma (1 + \alpha )}\int_a^b {\left| {f(x) + g(x)}
\right|^p(dx)^\alpha } } \right)^{1 / p} \\
 &\ge \left( {\frac{1}{\Gamma (1 + \alpha )}\int_a^b {\left| {f(x)}
\right|^p(dx)^\alpha } } \right)^{1 / p} + \left( {\frac{1}{\Gamma (1 +
\alpha )}\int_a^b {\left| {g(x)} \right|^q(dx)^\alpha } } \right)^{1 / q}. \\
\end{split}
 \tag{3.2}
 \label{3.2}
\end{equation*}
equalities holding if and only if $f(x) = \lambda g(x)$, where $\lambda $ is
a constant.
\end{thm}
\textbf{Proof}. Let $M = \frac{1}{\Gamma (1 + \alpha )}\int_a^b {\left| {f(x)}
\right|^p(dx)^\alpha } ,
N = \frac{1}{\Gamma (1 + \alpha )}\int_a^b {\left| {g(x)}
\right|^q(dx)^\alpha } $ and

\[
W = \left( {\frac{1}{\Gamma (1 + \alpha )}\int_a^b {\left| {f(x)}
\right|^p(dx)^\alpha } } \right)^{1 / p} + \left( {\frac{1}{\Gamma (1 +
\alpha )}\int_a^b {\left| {g(x)} \right|^q(dx)^\alpha } } \right)^{1 / q}.
\]

By H\"{o}lder inequality, in view of $0 < p < 1$, we have

\begin{equation*}
\begin{split}
 &W = \frac{1}{\Gamma (1 + \alpha )}\int_a^b {(\left| {f(x)} \right|^pM^{1 /
p - 1} + \left| {g(x)} \right|^pN^{1 / p - 1})(dx)^\alpha } \\
& \le \frac{1}{\Gamma (1 + \alpha )}\int_a^b {(\left| {f(x) + g(x)}
\right|^p(M^{1 / p} + N^{1 / p})^{1 - p})(dx)^\alpha } \\
 &= W^{1 - p}\frac{1}{\Gamma (1 + \alpha )}\int_a^b {(\left| {f(x) + g(x)}
\right|^p)(dx)^\alpha } .\\
 \end{split}
 \tag{3.3}
 \label{3.3}
\end{equation*}
By inequality (3.3), we arrive to reverse Minkowski's inequality and the
theorem is completely proved.
\begin{thm}\label{Corollary 3.3} Let $f_j (x) \ge 0$ , $j = 1,2,\ldots m$,

(1) for $p > 1$, we have
\begin{equation*}
\left( {\frac{1}{\Gamma (1 + \alpha )}\int_a^b {\left| {\sum\limits_{j =
1}^m {f_j (x)} } \right|^p(dx)^\alpha } } \right)^{1 / p} \le \sum\limits_{j
= 1}^m {\left( {\frac{1}{\Gamma (1 + \alpha )}\int_a^b {\left| {f(x)}
\right|^p(dx)^\alpha } } \right)^{1 / p}} .
\tag{3.4}
\label{3.4}
\end{equation*}

(2) for $0 < p < 1$ ,we have
\begin{equation*}
\left( {\frac{1}{\Gamma (1 + \alpha )}\int_a^b {\left| {\sum\limits_{j =
1}^m {f_j (x)} } \right|^p(dx)^\alpha } } \right)^{1 / p} \ge \sum\limits_{j
= 1}^m {\left( {\frac{1}{\Gamma (1 + \alpha )}\int_a^b {\left| {f(x)}
\right|^p(dx)^\alpha } } \right)^{1 / p}} .
\tag{3.5}
\label{3.5}
\end{equation*}
\end{thm}
\textbf{Proof}. (1) it follows from theorem 2.1 that

\begin{equation*}
\begin{split}
& \frac{1}{\Gamma (1 + \alpha )}\int_a^b {\left| {\sum\limits_{j = 1}^m {f_j
(x)} } \right|^p(dx)^\alpha } \le \frac{1}{\Gamma (1 + \alpha )}\int_a^b
{\vert f_1 (x)\vert \left| {\sum\limits_{j = 1}^m {f_j (x)} } \right|^{p -
1}(dx)^\alpha } \\
& + \cdots + \frac{1}{\Gamma (1 + \alpha )}\int_a^b {\vert f_m (x)\vert
\left| {\sum\limits_{j = 1}^m {f_j (x)} } \right|^{p - 1}(dx)^\alpha } \\
 &\le \left( {\frac{1}{\Gamma (1 + \alpha )}\int_a^b {\vert f_1 (x)\vert
^p(dx)^\alpha } } \right)^{1 / p}\left( {\frac{1}{\Gamma (1 + \alpha
)}\int_a^b {\left| {\sum\limits_{j = 1}^m {f_j (x)} } \right|^{q(p -
1)}(dx)^\alpha } } \right)^{1 / q} \\
& + \cdots + \left( {\frac{1}{\Gamma (1 + \alpha )}\int_a^b {\vert f_m
(x)\vert ^p(dx)^\alpha } } \right)^{1 / p}\left( {\frac{1}{\Gamma (1 +
\alpha )}\int_a^b {\left| {\sum\limits_{j = 1}^m {f_j (x)} } \right|^{q(p -
1)}(dx)^\alpha } } \right)^{1 / q} \\
 &= \sum\limits_{j = 1}^m {\left( {\frac{1}{\Gamma (1 + \alpha )}\int_a^b
{\vert f_j (x)\vert ^p(dx)^\alpha } } \right)^{1 / p}} \left(
{\frac{1}{\Gamma (1 + \alpha )}\int_a^b {\left| {\sum\limits_{j = 1}^m {f_j
(x)} } \right|^p(dx)^\alpha } } \right)^{1 / q} .\\
 \end{split}
 \tag{3.6}\label{3.6}
 \end{equation*}
in (3.6), multiplying both sides by $\left( {\frac{1}{\Gamma (1 + \alpha
)}\int_a^b {\left| {\sum\limits_{j = 1}^m {f_j (x)} } \right|^p(dx)^\alpha }
} \right)^{1 / q}$ yields (3.4).

(2) Similar to the proof of (3.2), we obtain (3.5).

\begin{cor}
\label{Corllary 3.4 }Let $f_j (x) \ge 0$，$j = 1,2,\ldots m$,

(1) for $p>1$, we have
\begin{equation*}
\frac{1}{\Gamma (1 + \alpha )}\int_a^b {\left| {\sum\limits_{j = 1}^m {f_j
(x)} } \right|^p(dx)^\alpha } > \sum\limits_{j = 1}^m {\frac{1}{\Gamma (1 +
\alpha )}\int_a^b {\left| {f_j (x)} \right|^p(dx)^\alpha } } .
\tag{3.7}
\label{3.7}
\end{equation*}

(2) for $0<p<1$, we have
\begin{equation*}
\frac{1}{\Gamma (1 + \alpha )}\int_a^b {\left| {\sum\limits_{j = 1}^m {f_j
(x)} } \right|^p(dx)^\alpha } < \sum\limits_{j = 1}^m {\frac{1}{\Gamma (1 +
\alpha )}\int_a^b {\left| {f_j (x)} \right|^p(dx)^\alpha } } \quad .
\tag{3.8}
\label{3.8}
\end{equation*}
\end{cor}
\begin{thm}
\label{Theorem 3.5} Let $f(x)$, $g(x) \ge 0$ , $0 < r < 1 < p$ , then
\begin{equation*}
\begin{split}
 &\left( {\frac{\frac{1}{\Gamma (1 + \alpha )}\int_a^b {\left| {f(x) + g(x)}
\right|^p(dx)^\alpha } }{\frac{1}{\Gamma (1 + \alpha )}\int_a^b {\left|
{f(x) + g(x)} \right|^r(dx)^\alpha } }} \right)^{1 / (p - r)} \\
 &\le \left( {\frac{\frac{1}{\Gamma (1 + \alpha )}\int_a^b {\left| {f(x)}
\right|^p(dx)^\alpha } }{\frac{1}{\Gamma (1 + \alpha )}\int_a^b {\left|
{f(x)} \right|^r(dx)^\alpha } }} \right)^{1 / (p - r)} + \left(
{\frac{\frac{1}{\Gamma (1 + \alpha )}\int_a^b {\left| {g(x)}
\right|^p(dx)^\alpha } }{\frac{1}{\Gamma (1 + \alpha )}\int_a^b {\left|
{g(x)} \right|^r(dx)^\alpha } }} \right)^{1 / (p - r)}.\\
 \end{split}
 \tag{3.9}
 \label{3.9}
\end{equation*}
equalities holding if and only if $f(x) = \lambda g(x)$.
\end{thm}
\textbf{Proof. }By Theorem 2.1 and Theorem 3.1, We have
\begin{equation*}
\begin{split}
 &  {\left( {\frac{1}
{{\Gamma (1 + \alpha )}}\int_a^b {{{\left| {f(x) + g(x)} \right|}^p}{{(dx)}^\alpha }} } \right)^{1/(p - r)}} \hfill \\
  & \leqslant {\left( {{{\left( {\frac{1}
{{\Gamma (1 + \alpha )}}\int_a^b {{{\left| {f(x)} \right|}^p}{{(dx)}^\alpha }} } \right)}^{1/p}} + {{\left( {\frac{1}
{{\Gamma (1 + \alpha )}}\int_a^b {{{\left| {g(x)} \right|}^p}{{(dx)}^\alpha }} } \right)}^{1/p}}} \right)^{p/(p - r)}} \hfill \\
   &= \left( {{{\left( {\frac{{\int_a^b {{{\left| {f(x)} \right|}^p}{{(dx)}^\alpha }} }}
{{\int_a^b {{{\left| {f(x)} \right|}^r}{{(dx)}^\alpha }} }}} \right)}^{1/p}}{{\left( {\frac{1}
{{\Gamma (1 + \alpha )}}\int_a^b {{{\left| {f(x)} \right|}^r}{{(dx)}^\alpha }} } \right)}^{1/p}}} \right. \hfill \\
 & {\left. { + {{\left( {\frac{{\int_a^b {{{\left| {g(x)} \right|}^p}{{(dx)}^\alpha }} }}
{{\int_a^b {{{\left| {g(x)} \right|}^r}{{(dx)}^\alpha }} }}} \right)}^{1/p}}{{\left( {\frac{1}
{{\Gamma (1 + \alpha )}}\int_a^b {{{\left| {g(x)} \right|}^r}{{(dx)}^\alpha }} } \right)}^{1/p}}} \right)^{p/(p - r)}} \hfill \\
  & \leqslant \left( {{{\left( {\frac{{\frac{1}
{{\Gamma (1 + \alpha )}}\int_a^b {{{\left| {f(x)} \right|}^p}{{(dx)}^\alpha }} }}
{{\frac{1}
{{\Gamma (1 + \alpha )}}\int_a^b {{{\left| {f(x)} \right|}^r}{{(dx)}^\alpha }} }}} \right)}^{1/(p - r)}} + {{\left( {\frac{{\frac{1}
{{\Gamma (1 + \alpha )}}\int_a^b {{{\left| {g(x)} \right|}^p}{{(dx)}^\alpha }} }}
{{\frac{1}
{{\Gamma (1 + \alpha )}}\int_a^b {{{\left| {g(x)} \right|}^r}{{(dx)}^\alpha }} }}} \right)}^{1/(p - r)}}} \right) \hfill \\
   &\times {\left( {{{\left( {\frac{1}
{{\Gamma (1 + \alpha )}}\int_a^b {{{\left| {f(x)} \right|}^r}{{(dx)}^\alpha }} } \right)}^{1/r}} + {{\left( {\frac{1}
{{\Gamma (1 + \alpha )}}\int_a^b {{{\left| {g(x)} \right|}^r}{{(dx)}^\alpha }} } \right)}^{1/r}}} \right)^{r/(p - r)}}.
\end{split}
 \tag{3.10}
 \label{3.10}
\end{equation*}
Using reverse Minkowski inequality implies that
\begin{equation*}
\begin{split}
 &\left( {\left( {\frac{1}{\Gamma (1 + \alpha )}\int_a^b {\left| {f(x)}
\right|^r(dx)^\alpha } } \right)^{1 / r} + \left( {\frac{1}{\Gamma (1 +
\alpha )}\int_a^b {\left| {g(x)} \right|^r(dx)^\alpha } } \right)^{1 / r}}
\right)^r \\
 &\le \frac{1}{\Gamma (1 + \alpha )}\int_a^b {\left| {f(x) + g(x)}
\right|^r(dx)^\alpha }. \\
 \end{split}
 \tag{3.11}
 \label{3.11}
\end{equation*}
By (3.10) and (3.11), we get (3.9). Hence, the theorem is completely proved.
\begin{cor}
\label{Corollary 3.6} Let $f_j (x) \ge 0$, $0 < r < 1 < p$ ,$j = 1,2,\ldots
m$ then
\begin{equation*}
\left( {\frac{\frac{1}{\Gamma (1 + \alpha )}\int_a^b {\left| {\sum\limits_{j
= 1}^m {f_j (x)} } \right|^p(dx)^\alpha } }{\frac{1}{\Gamma (1 + \alpha
)}\int_a^b {\left| {\sum\limits_{j = 1}^m {f_j (x)} } \right|^r(dx)^\alpha }
}} \right)^{1 / (p - r)} < \sum\limits_{j = 1}^m {\left(
{\frac{\frac{1}{\Gamma (1 + \alpha )}\int_a^b {\left| {f_j (x)}
\right|^p(dx)^\alpha } }{\frac{1}{\Gamma (1 + \alpha )}\int_a^b {\left| {f_j
(x)} \right|^r(dx)^\alpha } }} \right)^{1 / (p - r)}} .
\tag{3.12}
\label{3.12}
\end{equation*}
\end{cor}

 \end{document}